\documentclass{amsproc}% \usepackage{graphicx} \usepackage{amscd}
\usepackage{amsmath} \usepackage{amsfonts} \usepackage{amssymb}%
\theoremstyle{plain}

\numberwithin{equation}{section}

\begin{document}
\title[Isomorphisms between strongly triangular matrix rings]{Isomorphisms between strongly triangular matrix rings}
\author{P.N. \'Anh}
\address{R\'enyi Institute of Mathematics, Hungarian Academy of Sciences,
1364 Budapest, Pf.~127, Hungary}
\email{anh@renyi.hu}
\author{L. van Wyk}
\address{Department of Mathematical Sciences, Stellenbosch University,
P/Bag X1, \hfill\break Matieland~7602, Stellenbosch, South Africa }
\email{LvW@sun.ac.za}
\thanks{Corresponding author: L. van Wyk \\ \indent The research was carried
  out in accordance with the Hungary / South Africa Agreement on Cooperation in
  Science and Technology. In particular, the first author was
supported by the Hungarian National Foundation for Scientific Research under
Grants no. K-101515 and K-81928, and the second author was supported by the
National Research Foundation of South Africa under grant no.~UID 72375. Any
opinion, findings and conclusions or recommendations expressed in this
material are those of the authors and therefore the National Research
Foundation does not accept any liability in regard thereto. }
\thanks{The authors are very grateful to the anonymous referee of a previous version of this work. He/She pointed out some 
inaccuracies and called our attention to results in \cite{bhkp}. His/Her critical
remarks helped to improve this work. 
\\ \indent The authors also thank L. M\'arki and J. Szigeti for fruitful consultations.}

\begin{abstract} We describe isomorphisms between strongly triangular matrix rings that were defined earlier in 
\cite{bhkp} as ones having a complete set of triangulating idempotents, and we show that the so-called triangulating 
 idempotents behave analogously to idempotents in semiperfect rings. This
 study yields also a way to compute  theoretically the automorphism groups of
 such rings in terms of corresponding automorphism groups of certain subrings
 and bimodules involved in their structure, which completes the project started in \cite{avw}.

\end{abstract}

\maketitle

\noindent {\bf 1. Introduction}

\bigskip

\noindent Triangular matrix rings appear naturally in the theory of certain
algebras, like nilpotent and solvable Lie 
algebras, Kac-Moody, Virasoro and Heisenberg algebras (see, for example,
\cite{mp}), as well as in algebras of certain 
directed trees. In the latter case the triangular matrix rings may be seen to
provide the abstract description of such 
quiver algebras without mentioning the associated directed tree and without appropriate numbering of the vertices.

Triangular matrix rings have become an important object of intense research,
for example, it is a key tool in the 
description of semiprimary hereditary rings (see, for example, \cite{f}), and
certain triangular matrix rings are 
natural examples of representation-finite hereditary algebras (see, for example, \cite{begp} and \cite{g}). 

On the other hand, Birkenmeier et al in \cite{bhkp} developed the general
theory of generalized triangular matrix rings 
and used it to describe several particular classes of rings. Combining their
terminology with ones (introduced later) 
in \cite{avw} we say that a ring $A$ admits an 
\emph{$m$-strongly (upper) triangular matrix decomposition with respect to an
  ordered sequence} $\{e_1,\ldots,e_m\}$ 
if the $e_i$'s are  pairwise orthogonal idempotents in $A$ such that
$1=e_1+\cdots +e_m, \ e_jAe_i=0$ for all $j>i$ and $e_iAe_i$ is semicentral
reduced for every $i$, or equivalently, $\{e_1,\ldots,e_m\}$ is 
\emph{a complete set of left triangulating idempotents} by terminology of
\cite{bhkp} . Here, according to \cite{bhkp}, an idempotent~$e$ in a
ring~$A$ is called {\it semicentral} if $(1-e)Ae=0,$ and $A$ is called {\it
  semicentral reduced} if 
$0$ and~$1$ are the only semicentral idempotents in
$A$, i.e., $A$ is semicentral reduced if and only if $A$ is strongly
indecomposabble in the sense of~\cite{avw}. Therefore, an idempotent $e\in A$
is semicentral reduced if it is semicentral and the subring $eAe$ is a
strongly indecomposable ring. If we set 
$R_i := e_i A e_i$ and
$L_{ij} := e_iAe_j$ for $i < j$, then $A$ can be written as a generalized
upper triangular matrix ring 

\smallskip

$$\left[ \begin{array}{ccccc}
R_1 & L_{12} & L_{13} & \cdots & L_{1m} \\
\\
0 & R_2 & L_{23} & \cdots & L_{2m} \\
\\
\vdots & \ddots & \ddots &  & \vdots \\
\\
0 &  \cdots & 0 & R_{m-1} & L_{m-1, m} \\
\\
0 & \cdots & \cdots & 0 & R_m \end{array} \right]$$ 

\medskip

\noindent with the obvious matrix addition and multiplication. 
 It was pointed out
in \cite{bhkp} that by reversing the 
order of the sequence $\{e_1,\ldots,e_m\}$ one obtains a new sequence
providing the 
lower triangular matrix representation for $A$. Therefore it is not a
restriction to study rings with a complete set of left triangulating idempotents.

The aim of this paper is to describe isomorphisms between strongly triangular matrix rings, thereby finishing the project initiated in \cite{avw}. As a by-product we show that triangulating idempotents behave similarly to idempotents in semiperfect rings. Namely, if one fixes a complete set $\{e_1,\ldots,e_m\}$
 of triangulating idempotents, then a left ideal generated by any semicentral
 idempotent is isomorphic to one generated by 
an appropriate partial sum of some idempotents from the set $\{e_1,\ldots,e_m\}$.

For more information and detailed treatment of triangular matrix rings and their applications in other areas of mathematics we refer to \cite{bhkp}, and for some interesting related questions on matrix rings we refer to \cite{s}.
 
\bigskip

\noindent {\bf 2.   Strongly triangular matrix rings}

\bigskip

\noindent A strongly (upper) triangular matrix decomposition of a ring $A$
depends essentially on the ordered sequence $\{e_1, \ldots, e_m\}$ of pairwise
orthogonal idempotents with sum~1. However, in particular cases, another
ordering of the set $\{e_1, \ldots, e_m\}$ may also give a strongly triangular
matrix decomposition of $A$. 

Furthermore, if there is no room for
misunderstanding, then for short we sometimes say that a
ring $A$ is a {\it strongly triangular matrix ring}, without stating exactly
the ordering on the set $\{e_1, \ldots, e_m\}$. Therefore one has to see clearly that all $R_i$ are
semicentral reduced, but all $e_i$, except $e_1$, need not be even semicentral
idempotents of $A$, i.e., $e_i$ for $i\geqq 1$ is certainly reduced semicentral
only in the subring 
$A_i=(e_i+\cdots +e_m)A(e_i+\cdots +e_m)$ of $A=A_1$ but not necessarily in $A_j$ with $j<i$. For example, if $A$ is a
strongly triangular matrix ring with respect to the ordered sequence $\{e_1,
e_2, e_3\}$, then the generalized matrix decompositions of $A$ with respect 
to the
ordered sequence $\{e_2, e_1, e_3\}$ and $\{e_2, e_3, e_1\}$ are

$$\left[ \begin{array}{ccccc}
R_2 & 0 & L_{23} \\
\\
L_{12} & R_1 & L_{13}\\
\\
0 & 0 & R_3 \end{array} \right]\qquad
{\rm{\text and}}\qquad \left[ \begin{array}{ccccc}
R_2 & L_{23} & 0 \\
\\
0 & R_3 & 0\\
\\
L_{12} & L_{13} & R_1 \end{array} \right],$$

\bigskip

\noindent respectively, which are definitely not triangular matrix
decompositions of~$A$.

Next, let $B$ be an $n$-strongly triangulated matrix ring with respect to an ordered sequence $\{f_1, \ldots, f_n\}$, i.e., the $f_i$'s are pairwise orthogonal idempotents in $B$ with sum 1, $f_jBf_i=0$ for all $j>i, \ S_i := f_iBf_i$ is semicentral reduced for every $i$, and $f_i$  is a semicentral reduced idempotent of the ring $B_i = (f_i + \cdots + f_n) B(f_i+ \cdots + f_n)$ for $i=1,\ldots, n$. For each $i \neq j$, let $M_{ij} = f_i Bf_j$. Therefore $M_{ij} =0$ for all $j < i$. Moreover, for each $i$ let $M_i$ be the $i$-th truncated row of $B$, i.e.,
$$M_i = \displaystyle\oplus_{k>i} M_{ik} = \displaystyle\oplus_{k\neq i} M_{ik}.$$
If $\sigma$ is any permutation on $\{1, \ldots, n\}$, then $\sigma$ induces a new (generalized) matrix ring decomposition on $B$ with respect to the ordered sequence $\{f^\sigma_1:=f_{\sigma(1)}, \ldots, f^\sigma_n:=f_{\sigma (n)}  \}$. According to this notation, if we write $g_i = f^\sigma_i$, then one can identify the above convention as follows. Let $T_i = g_iBg_i,  \ C_i = (g_i+ \cdots + g_n)B(g_i + \cdots + g_n), \ N_{ij} = g_iBg_j$ for all $i \neq j, \ N_i = \displaystyle\oplus_{k \neq i} N_{ik}$. It is important to emphasize that $B$ is not necessarily an $n$-strongly triangular matrix ring with respect to the ordered sequence $\{f^\sigma_1, \ldots, f^\sigma_n\}$. From the above one gets
$$T_i = S^{\sigma}_i:=S_{\sigma(i)}, \ C_i =  B^\sigma_i:=B_{\sigma (i)}, \ N_i =M^\sigma_i:= M_{\sigma (i)}.$$

\bigskip

Now we are in a position to state the main result precisely.

\bigskip

\noindent {\bf Theorem.} Let $A$ and $B$ be $m$- and $n$-strongly triangular
matrix rings with respect to ordered sequences $\{e_1, \ldots, e_m\}$ and
$\{f_1, \ldots, f_n\}$, respectively. Then $A$ and $B$ are isomorphic via an
isomorphism $\varphi: A \rightarrow B$ iff $m=n$ and there is a permutation
$\sigma$ of $\{1, \ldots, m\}$ such that $B$ is also an $m$-strongly
triangular matrix ring with respect to the ordered sequence $\{f^\sigma_1,
\ldots, f^\sigma_m\}$, there are ring isomorphisms $\rho_i: R_i \rightarrow S^\sigma_i = 
S_{\sigma (i)}, \ i=1,\ldots, m=n,$ and for $i = 1, \ldots, m-1$ there are
elements 
$m_i \in M^\sigma_i$ and ring isomorphisms $\varphi_{i+1}: A_{i+1} \rightarrow
B^\sigma_{i+1}$ and 
$R_i - A_{i+1}$-bimodule isomorphisms
$\chi_i : e_iA_i(e_{i+1} + \cdots + e_{m(=n)}) = L_i \rightarrow M^\sigma_i$
with respect to $\rho_i, \varphi_{i+1}$, such that 
for $i =1, \ldots, m-1$ and 

$$a_i = \left[ \begin{array}{rrr} r_i & \ell_i \\
\\
0 & a_{i+1} \end{array} \right] \in A_i = \left[ \begin{array}{rrr}
R_i & L_i \\
\\
0 & A_{i+1} \end{array}\right],$$

\medskip 

$$\varphi_i(a_i) = \left[ \begin{array}{ccc}
\rho_i(r_i) & \rho_i(r_i)m_i+ \chi_i (\ell_i) - m_i \varphi_{i+1}(a_{i+1})\\
\\
0 & \varphi_{i+1}(a_{i+1}) \end{array}\right].$$

\bigskip 

\noindent Moreover, all isomorphisms between isomorphic rings $A$ and $B$ can
be described in this manner. (Keep in mind that $\varphi_1 = \varphi, \ 
\varphi_m = \rho_m; \ A_m  = R_m.$)

\bigskip
\noindent{\bf Remark 1.} The equality $m=n$ as well as some invariants
associated to a complete set of triangulating idempotents up to a permutation $\sigma$ were already 
obtained as~Theorem~2.10 in~\cite{bhkp}, where an isomorphism is (surprisingly
enough) an inner automorphism. However, 
these results are by-products of our description of general isomorphisms
between such rings and our treatment is 
both elementary and direct. For further details for structural discussion we
refer to Theorems 2.10, 3.3 and 
Corollary 3.4 in~\cite{bhkp}.
 
\bigskip
\noindent {\bf Proof of Theorem.}
We prove this theorem by induction on $m$, a number of pairwise orthogonal
idempotents $e_i$ in the ordered sequence $\{e_1, \ldots, e_m\}$ giving a
strongly triangular matrix ring decomposition on $A$. The case $m=1$ is 
obvious by the definition, because $B$ must be also semicentral reduced, i.e., 
$m=n=1$. Assume now that $m\geq 2$ and the theorem holds for $m-1$.

The first induction step is the following obvious but interesting result 
(by direct computation, see also \cite{avw}). Because of its importance we
state it separately as a self-contained assertion.

\bigskip

\noindent {\bf Proposition.} Let $e \in A$ and $ f\in B$ be semicentral 
idempotents. Put $R = eAe, \ S = fBf, \bar A = (1-e)A(1-e), \bar B=(1-f)B(1-f)$, \ $L = eA(1-e), \ M=fB(1-f)$, i.e., $A = \left[ \begin{array}{ll} R& L \\
0 & \bar A \end{array} \right], \ B = \left[ \begin{array}{ll} S& M\\
0 & \bar B \end{array}\right]$, and let $\varphi : A \rightarrow B$ be a ring isomorphism. Then $\varphi (e) \in f+M$ if and only if there are ring isomorphisms $\rho: R \rightarrow S$ and $\bar{\varphi}: \bar A \rightarrow \bar B$ and an $R - \bar A$-bimodule isomorphism $\chi : L \rightarrow M$ ($M$ is an $R-\bar A$-bimodule via $\rho$ and $\bar{\varphi}$) and an element $m$ in $M$ such that

\medskip

\hspace*{1cm}$\varphi \left(\left[ \begin{array}{lll} r & \ell \\
\\
0 & a \end{array}\right]\right) = \left[ \begin{array}{ccc}
\rho (r) & \rho (r) m + \chi (\ell ) - m \bar{\varphi} (a)\\
\\
0 & \bar{\varphi} (a) \end{array}\right].\hfill (1)$

\medskip

\noindent In particular, $\chi$ is just the restriction of $\varphi$ to $L$. Moreover, all isomorphisms $\varphi$ from~$A$ to $B$ satisfying $\varphi (e)\in f+M$ can be obtained from a quadruple $(\rho,\bar{\varphi},\chi,m)$ in this manner.

\medskip

For the verification of the {\bf Proposition} one observes $\varphi (e)\in f+M$ if and only if 
$\varphi (e)=f+m=f_m$ for some $m\in M$. Therefore $\varphi (1-e)=1-f-m=
(1-f)-m=g_m$. Put $g=1-f$. Then by direct calculations (see also Lemma 2.2 in \cite{avw}) one has $M=fBg=f_mBg_m$
and canonical isomorphisms 
$S\cong f_mBf_m: v\in S \mapsto v+vm\in f_mBf_m$ and $\bar B\cong g_mBg_m: 
w\in B_1 \mapsto w-mw\in g_mBg_m$. Consequently,
$\varphi$ induces the isomorphisms $\rho: R\longrightarrow S:r\in R\mapsto
\rho(r)=v, \ 
\bar{\varphi}:\bar A\longrightarrow \bar B:a\in \bar A\mapsto \bar{\varphi}(a)=w$
if $\varphi(r)=v+vm\in f_mBf_m=\varphi(e)\varphi(A)\varphi(e)=f_mBf_m, \varphi(a)=w-mw\in 
\varphi(1-e)\varphi(A)\varphi(1-e)=g_mBg_m$. Therefore for an arbitrary element of $A$,
i.e., for an arbitrary triple $r\in R, \ell \in L, a\in A$, one obtains~(1) immediately. 

Finally, it is clear that every quadruple $(\rho,\bar{\varphi},\chi,m)$ as described in the statement of the proposition leads to one of the desired isomorphisms, 
completing the justification of the proposition. 

\bigskip

Now we continue with the proof of the {\bf Theorem}. Consider the

\bigskip

\noindent {\bf Main Step.} Let $A$ and $B$ be $m$-and $n$-strongly upper triangular matrix rings with respect to $\{e_1, \ldots, e_m\} \subseteq A$ and $\{f_1, \ldots, f_n\} \subseteq B$, respectively, and let $\varphi: A \rightarrow B$ be a ring isomorphism. Let $R_i = e_i Ae_i, \ L_{ij} =e_i Ae_j$ for $i < j$, and $S_i = f_iBf_i, \ M_{ij} = f_i Bf_j$ for $i < j$, i.e.,
\medskip

\hspace*{-0.5cm}$
A = \left[ \begin{array}{ccccc}
R_1 & L_{12} & L_{13} & \cdots & L_{1m} \\
\\
0 & R_2 & L_{23} & \cdots & L_{2m} \\
\\
\vdots & \ddots & \ddots &  & \vdots \\
\\
\vdots &  & \ddots & \ddots &  \vdots \\
\\
0 & \cdots & \cdots & 0 & R_m \end{array} \right], \ \ 
B= \left[ \begin{array}{ccccc}
S_1 & M_{12} & M_{13} & \cdots & M_{1n} \\
\\
0 & S_2 & M_{23} & \cdots & M_{2n} \\
\\
\vdots & \ddots & \ddots &  & \vdots \\
\\
\vdots &  & \ddots & \ddots & \vdots \\
\\
0 & \cdots & \cdots & 0 & S_n \end{array} \right].\hfill$

\medskip

\noindent Then either $\varphi (e_1) \in f_1 + M_1$ or there is a $j \geq 2$
such that $\varphi (e_1) \in f_j + M_j$ and $M_{1j} =0, \ldots, M_{j-1,j} =0$.

\bigskip

Since $f = \varphi (e_1) \in B$ is a semicentral reduced idempotent, the statement of the {\bf Main Step} can be reformulated in an equivalent, but little sharper, form, namely: 

\bigskip

If $f$ is a semicentral reduced idempotent in the $n$-strongly triangulated matrix ring~$B$, then either $f \in f_1 + M_1$ or there is a $j \geq 2$ such that $f\in f_j + M_j$ and $M_{1j} =0, \ldots, M_{j-1, j} =0$.

\bigskip 

\noindent {\bf Proof of the Main step.} Again we use induction for the
verification. Let $F_1 = 1-f_1, \ B_2 = F_1 BF_1, \ M = f_1BF_1$, i.e., 
$B = \left[ \begin{array}{ll} S_1 & M\\
0 & B_2 \end{array} \right]$. The statement is obvious for $n=1$. Assume $n \geq 2$.  Writing $f = \left[ \begin{array}{lll} \alpha & \mu \\
0 & \beta \end{array} \right] \in B = \left[ \begin{array}{lll} S_1 & M \\
0 & B_2 \end{array}\right]$ it follows from $f = f^2 = \left[ \begin{array}{ccc} \alpha^2 & \alpha \mu + \mu \beta\\ \\ 0 & \beta^2 \end{array}\right]$ that $\alpha^2 = \alpha, \ \beta^2 = \beta$ and $\alpha \mu + \mu \beta = \mu$, and so $\alpha \mu \beta =0$. Writing $s = \left[ \begin{array}{lll} \alpha & \alpha\mu \\
0 & 0 \end{array}\right]$ and $b = \left[ \begin{array}{lll} 0 & \mu\beta\\
0 & \beta \end{array}\right]$ we get $s^2 =s, \ b^2 = b$ and $sb = 0 = bs$. Hence, $f = s+b$ implies that $fs = s = sf$ and $fb = b = fb$, i.e., $s,b \in fBf$, which is semicentral reduced.  Moreover,
$$\left[ \begin{array}{ll} 0 & \mu\beta \\
0 & \beta \end{array}\right] \ \left[ \begin{array}{ll} S_1 & M\\
0 & B_2 \end{array}\right] \ \left[ \begin{array}{ll} \alpha & \alpha\mu\\
0 & 0 \end{array}\right] = \left[ \begin{array}{ll} 0 & \mu\beta B_2 \\
0 & \beta B_2 \end{array}\right] \ \left[ \begin{array}{ll} \alpha & \alpha\mu\\
\alpha & 0 \end{array}\right] =0,$$
i.e., $s \in f B f$ is a semicentral. Consequently,
$$f =s= \left[ \begin{array}{ll} \alpha & \alpha\mu \\
0 & 0 \end{array}\right] \quad \mbox{or} \quad f = \left[ \begin{array}{ll} 0 & \mu\beta \\
0 & \beta \end{array}\right].$$

Assume the first case: $f = \left[ \begin{array}{ll} \alpha & \mu \\
0 & 0 \end{array}\right]$. Then $(1-f)Bf =0$ implies that
$$\left[ \begin{array}{ccc} 1-\alpha & -\mu \\
0 &1 \end{array} \right] \ \left[ \begin{array}{ll} S_1 & M \\
0 & B_2 \end{array}\right] \ \left[ \begin{array}{ll} \alpha & \mu \\
0 & 0 \end{array} \right] = \left[ \begin{array}{ccc} (1-\alpha )S_1\alpha & \ast \\
0 & 0 \end{array} \right]=0,$$
i.e., $(1-\alpha) S_1 \alpha =0$, hence $\alpha$ is a semicentral idempotent in $S_1$. Since $\alpha \neq 0$ and $S_1$ is semicentral reduced we obtain $\alpha =1$, i.e., $f \in f_1B$. 

Next, consider the case $f = \left[ \begin{array}{ll} 0 & \mu \\
0 & \beta \end{array}\right]$. Again $(1-f)Bf =0$ implies that
$$\left[ \begin{array}{cc} 1 & -\mu \\
0 & 1- \beta \end{array}\right] \left[ \begin{array}{cc} S_1 & M \\
0 & B_2 \end{array}\right] \left[ \begin{array}{ll} 0 & \mu \\
0 & \beta \end{array}\right] = \left[ \begin{array}{cc} 0 & \ast \\
0 & (1-\beta)B_2 \beta \end{array}\right] =0,$$
showing that $(1-\beta ) B_2 \beta=0$, i.e., $\beta$ is semicentral. Since $\beta B_2 \beta = \beta B \beta = f B f$, we have that $\beta$ is also reduced. Since $B_2$ is $(n-1)$-strongly triangular, the induction hypothesis shows that there is  a $j \geq 2$ such that, in 
$B_2 = \left[ \begin{array}{cccc} S_2 & M_{23} & \cdots & M_{2n} \\
0 & S_3 & \cdots & M_{3n}\\
 \vdots& \ddots & \ddots & \vdots\\
0& \cdots & 0 & S_n \end{array}\right]$, $f$ is of the form
$$j-1 \left[ \begin{array}{ccccccc} 
 0 & \cdots & 0 & 0 &&&  \\
 & \ddots & \vdots & \vdots & & \bigcirc & \\
 & & 0 & 0 &&&  \\
&& & 1  & \ast & \cdots & \ast \\
  &&&& 0 & \cdots & 0  \\
  &&&&  &  \ddots &  \vdots \\
   &&&&&& 0 \\
   \end{array}\right],$$ i.e., $M_{2j} =0, \ldots, M_{j-1,j} =0.$
 Therefore, in $B$,  \ $f$ is of the form
 $$j \left[ \begin{array}{ccccccccc} 
 0 & \cdots & \cdots & 0 & x_1 &&&&  \\
 & \ddots && \vdots &0&  & & & \\
 & & \ddots & \vdots & \vdots  &  & & \bigcirc &  \\
 &&& 0 & 0 &&&& \\
 &&& & 1  & \ast & \cdots & \cdots & \ast \\
 &&&&&  0 & \cdots & \cdots & 0  \\
  &&&&  & &  \ddots & & \vdots \\
  &&&&  & & & \ddots & \vdots \\
   &&&&&&&& 0 \\
   \end{array}\right],$$ where $x_1 \in M_{1j}$.
  Now $0=(1-f)Bf=$ 
  
 \medskip
  
 \hspace*{-1.2cm}$ \left[ \begin{array}{ccccccccc} 
 1 & 0 & \cdots & 0 & -x_1 &&&&  \\
 & \ddots & &  &0&  & & & \\
 & & \ddots &  & \vdots  &  & & \bigcirc &  \\
 &&& 1 & 0 &&&& \\
 &&& & 0  & \ast & \cdots & \cdots & \ast \\
 &&&&&  1 & 0 & \cdots & 0  \\
  &&&&  & &  \ddots & \ddots & \vdots \\
  &&&&  & & & \ddots & 0 \\
   &&&&&&&& 1 \\
   \end{array}\right] B
    \left[ \begin{array}{ccccccccc} 
 0 & \cdots & \cdots & 0 & x_1 &&&&  \\
 & \ddots && \vdots &0&  & & & \\
 & & \ddots & \vdots & \vdots  &  & & \bigcirc &  \\
 &&& 0 & 0 &&&& \\
 &&& & 1  & \ast & \cdots & \cdots & \ast \\
 &&&&&  0 & \cdots & \cdots & 0  \\
  &&&&  & &  \ddots & & \vdots \\
  &&&&  & & & \ddots & \vdots \\
   &&&&&&&& 0 \\
   \end{array}\right]\hfill$ 
  
\medskip
 
\noindent implies both $x_1 =0$ and $M_{1j} =0$, completing the proof of
the {\bf Main Step}. \ $\square$

\medskip

The following observation is the last piece in the proof 
of the {\bf Theorem}.

\medskip

\noindent {\bf Lemma.} If $$B= \left[ \begin{array}{ccccc}
S_1 & M_{12} & M_{13} & \cdots & M_{1n} \\
\\
0 & S_2 & M_{23} & \cdots & M_{2n} \\
\\
\vdots & \ddots & \ddots &  & \vdots \\
\\
\vdots &  & \ddots & \ddots & \vdots \\
\\
0 & \cdots & \cdots & 0 & S_n \end{array} \right]$$ 
is $n$-strongly triangular with respect to the ordered sequence 
$\{f_1, \cdots, f_n\}$ of pairwise orthogonal idempotents such that  
$M_{1j} =0, \ldots, M_{j-1,j} =0$  for some index $j > 1$, then $B$ is 
also $n$-strongly triangular with respect to the ordered sequence 
$\{f_j, f_1, \ldots, f_{j-1}, f_{j+1}, \ldots, f_n\}$.

\medskip

\noindent {\bf Proof.} Obvious by definition. \ $\square$

\medskip

If we define now $\sigma(1)=j$, then the above Lemma together with the 
Proposition
shows that $\varphi$ induces the ring isomorphisms 
$\rho_1: R_1\cong S_j=S^\sigma_1 = S_{\sigma (1)}, \ \varphi_2:A_2\cong
\bar B=(1-f_j)B(1-f_j)$ and the bimodule isomorphism 
$\chi_1:L_1=e_1A(1-e_1)\cong M^\sigma_1=M_{\sigma (1)}=f_jB(1-f_j)$
together with an element $m_1\in M^\sigma_1$ such that for an arbitrary
$a=a_1 = \left[ \begin{array}{rrr} r_1 & \ell_1 \\
0 & a_2 \end{array} \right] \in A=A_1 = \left[ \begin{array}{rrr}
R_1 & L_1 \\
0 & A_2 \end{array}\right],$ $\varphi=\varphi_1$ 
satisfies

$$\varphi(a)=\varphi_1(a_1) = \left[ \begin{array}{ccc}
\rho_1(r_1) & \rho_1(r_1)m_1+ \chi_1 (\ell_1) - m_1 \varphi_2(a_2)\\
\\
0 & \varphi_2(a_2) \end{array}\right],$$

\noindent and every such $\varphi$ can be described in this manner. Since $A_1$ is an $(m-1)$-strongly triangular matrix ring and $\bar B$
is an $(n-1)$-strongly triangular matrix ring, the theorem follows now
immediately from the induction hypothesis which makes the proof of the {\bf Theorem} complete. \ $\square$

\medskip

We  emphasize  three important remarks.

\medskip

\noindent {\bf Remark 2.} If $M_{i,i+1} \neq 0$ for $i =1, \ldots, n-1$, then $\{f_1, \ldots, f_n\}$ is the unique order (up to isomorphism) which induces the $n$-strongly triangular matrix decomposition on $B$.

\medskip

\noindent {\bf Remark 3.} If $i<j$ and $M_{ij} =0, \ldots,M_{j-1,j}=0,$ then  $\{f_1, \ldots, f_{i-1}, f_j, f_i, f_{i+1}, \ldots,\hfill\break f_{j-1}, f_{j+1}, \ldots, f_n\}$ also induces an $n$-strongly triangular matrix decomposition on~$B$. In this case, the ring $(f_i+\cdots +f_{j-1}+f_j)B(f_i+\cdots +f_{j-1}+f_j)$ is the direct sum of the two rings $(f_i+\cdots +f_{j-1})B(f_i+\cdots +f_{j-1})$ and $S_j=f_jBf_j$. In particular, in the case $i=1, j=n$ the ring $B$ is the direct sum of $(f_1+\cdots +f_{n-1})B(f_1+\cdots +f_{n-1})$ and~$S_n$ if the truncated last column is 0. 

\medskip

\noindent {\bf Remark 4.} Specializing the theorem for the case $A=B$ one
obtains the description of the automorphism group of the strongly triangular
matrix rings in terms of the corresponding automorphism groups of reduced
rings $R_i$ and of the corresponding bimodules $L_i$ similar to one given in
\cite{avw}.

\medskip
The proof of the main step shows also that for an arbitrary semicentral
idempotent $e$ in a strongly triangular 
matrix ring $A$ each semicentral reduced
idempotent $g\in A$ is either  $g = \left[ \begin{array}{ll} \alpha & \mu \\
0 & 0 \end{array}\right]$ or  $g = \left[ \begin{array}{ll} 0 & \nu \\
0 & \beta \end{array}\right]$ where $\alpha \in eAe$ and $\beta \in (1-e)A(1-e)$ are 
semicentral reduced idempotents in $A$, the associated subrings $gAg, \alpha A\alpha, \beta A\beta$ 
are isomorphic, and 
$\mu,\, \nu$ are appropriate elements in $eA(1-e)$. Observing that
$C=(1-\alpha)A(1-\alpha)$ in the first case or $C=(1-\beta)A(1-\beta)$ in the
second case is an $(m-1)$-strongly triangular matrix ring, by considering $\bar e=e-\alpha$ in the first case or in view of 
$e=(1-\beta)e(1-\beta)$ in the second case, respectively, the Main Step and the Theorem
together with an obvious induction imply immediately the following

\bigskip

\noindent {\bf Corollary.} Any semicentral idempotent $e$ in a $m$-strongly
triangular matrix ring~$A$ with a complete 
set of triangulating idempotents 
can be written as a sum of $l$ pairwise orthogonal idempotents
$\{e_1,\ldots,e_l\}$ where $l\leq m$ is uniquely determined by $e$ and this
set of idempotents can be extended to the 
first $l$ idempotents in a complete set of triangulating idempotents of $A$.

\end{document}